\documentclass[11pt]{amsart}
\usepackage{mathrsfs}
\usepackage[colorlinks,linkcolor=blue,hyperindex]{hyperref}

\numberwithin{equation}{section}

\usepackage{amsfonts}

\newcommand{\Om}{\Omega}
\newcommand{\bd}{\begin{enumerate}}
\newcommand{\ed}{\end{enumerate}}
\newcommand{\btheorem}{\begin{theorem}}
\newcommand{\etheorem}{\end{theorem}}
\newcommand{\bproposition}{\begin{proposition}}
\newcommand{\eproposition}{\end{proposition}}
\newcommand{\bdefinition}{\begin{definition}}
\newcommand{\edefinition}{\end{definition}}
\newcommand{\bcorollary}{\begin{corollary}}
\newcommand{\ecorollary}{\end{corollary}}
\newcommand{\bproof}{\begin{proof}}
\newcommand{\eproof}{\end{proof}}
\newcommand{\bremark}{\begin{remark}}
\newcommand{\eremark}{\end{remark}}
\newcommand{\eexample}{\end{example}}
\newcommand{\bexample}{\begin{example}}
\newcommand{\elemma}{\end{lemma}}
\newcommand{\blemma}{\begin{lemma}}
\newcommand{\ee}{\end{eqnarray*}}
\newcommand{\be}{\begin{eqnarray*}}
\newcommand{\p}{\partial}
\renewcommand{\bar}{\overline}

\begin{document}
\newtheorem{claim}{Claim}
\newtheorem{theorem}{Theorem}[section]
\newtheorem{lemma}[theorem]{Lemma}
\newtheorem{corollary}[theorem]{Corollary}
\newtheorem{definition}[theorem]{Definition}
\newtheorem{proposition}[theorem]{Proposition}
\newtheorem{question}{question}[section]
\newtheorem{example}[theorem]{Example}
\newtheorem{conjecture}[theorem]{Conjecture}

\theoremstyle{definition}
\newtheorem{remark}[theorem]{Remark}

\newcommand{\sq}{\sqrt{-1}}
\renewcommand{\phi}{\varphi}
\newcommand{\bp}{\bar{\partial}}
\newcommand{\beq}{\begin{equation}}
\newcommand{\eeq}{\end{equation}}

\newcommand{\G}{{\mathbb G}}
\renewcommand{\H}{{\mathbb H}}
\renewcommand{\P}{{\mathbb P}}

\renewcommand{\tilde}{\widetilde}
\newcommand{\ts}{\otimes}
\newcommand{\sX}{\mathcal X}
\newcommand{\sO}{\mathcal O}
\newcommand{\sE}{\mathcal E}
\newcommand{\sL}{\mathcal L}
\renewcommand{\>}{\rightarrow}
\renewcommand{\hat}{\widehat}

\title[Curvatures of direct image  of vector bundles and applications]{Curvatures of direct image sheaves of vector bundles and applications$^*$}
\author[Kefeng Liu]{Kefeng Liu$^\dagger$}

\address{$^\dagger$Department of Mathematics, UCLA, 520 Portola Plaza,
 Los Angeles, CA 90095}
\author[Xiaokui Yang]{Xiaokui Yang$^\ddagger$}
\thanks{$^{*}$Research support}
\address{$^\ddagger$Department of Mathematics, Northwestern University, 2033 Sheridan Road, Evanston, IL 60208}
\maketitle

\begin{abstract}  Let $p:\sX\>S$ be a proper K\"ahler fibration and $\sE\>\sX$ a Hermitian holomorphic vector bundle.
 As motivated by the work of Berndtsson(\cite{Berndtsson09a}), by using basic Hodge theory, we derive several general
curvature formulas for the direct image $p_*(K_{\sX/S}\ts \sE)$ for
general Hermitian holomorphic vector bundle $\sE$ in a  simple way.
 A straightforward application is that, if the family $\sX\>S$ is
 infinitesimally trivial and
Hermitian vector bundle $\sE$ is Nakano-negative along the
 base $S$, then the direct image $p_*(K_{\sX/S}\ts
\sE)$ is Nakano-negative. We also use these curvature formulas to
study the moduli space of projectively flat vector bundles with
positive first Chern classes and obtain that, if the Chern curvature
of direct image $p_*(K_{X}\ts E)$--of a positive projectively flat
family $(E,h(t))_{t\in \mathbb D}\>X$--vanishes, then the curvature
forms of this family are connected by holomorphic automorphisms of
the pair $(X,E)$.

\end{abstract}
%
%

\setcounter{tocdepth}{1} \tableofcontents

\section{Introduction}

Let $\sX$ be a K\"ahler manifold with dimension $m+n$ and $S$ a
K\"ahler manifold with  dimension $m$. Let $p:\sX\>S$  be a proper
K\"ahler fibration. Hence, for each $s\in S$,
$$X_s:=p^{-1}(\{s\})$$
is a compact K\"ahler manifold with dimension $n$. Let
$(\sE,h^{\sE})\>\sX$ be a Hermitian holomorphic vector bundle.
Consider the space of holomorphic $\sE$-valued $(n,0)$-forms on
$X_s$,
 $$E_s:=H^0(X_s,\sE_s\ts K_{X_s})\cong H^{n,0}(X_s,\sE_s)$$
where $\sE_s=\sE|_{X_s}$. It is well-known that, if the vector
bundle $E$ is  ``positive" in certain sense, there is a natural
holomorphic structure on
$$E=\bigcup_{s\in S}\{s\}\times E_s$$  such that $E$ is isomorphic to the direct
image sheaf $p_*(K_{\sX/S}\ts \sE)$.
 Using the canonical isomorphism
$K_{\sX/S}|_{X_s}\cong K_{X_s},$
 a \emph{local smooth} section $u$ of $E$ over $S$ can be identified as  a holomorphic
$\sE$-valued $(n,0)$ form on $X_s$.

By the identification above, there is a natural metric on $E$. For
any local smooth section $u$ of $E$, one can define a Hermitian
metric on $E$ by \beq h(u,u)=c_n\int_{X_s} \left\{u,u\right\}\eeq
where $c_n=(\sq)^{n^2}$. Here, we only use the Hermitian metric of
$\sE_s$ on each fiber $X_s$ and we do not specify  background
K\"ahler metrics on the fibers.  Berndtsson defined in \cite[Lemma
~4.1]{Berndtsson09a},  a natural Chern connection $D$ on $(E,h)$,
and  computed the curvature tensor of direct image $p_*(K_{\sX/S}\ts
\sL)$ of (semi-)positive line bundle $\sL\>\sX$.

Next we would like to describe our results in this paper briefly. As
motivated by the work of Berndtsson(\cite{Berndtsson09a}), we
compute the curvature tensor of the direct images $p_*(K_{\sX/S}\ts
\sE)$ for arbitrary Hermitian vector bundles $\sE\>\sX$ by using
basic Hodge theory which also simplify Berndtsson's original proofs
mildly. In the following formulations, if not otherwise stated, we
do not make any positivity or negativity assumption on the curvature
tensors of $\sE$ or $\sL$, but we assume that every $E_s$ has the
same dimension.

   Let $(X,\omega_g)$ be a compact K\"ahler manifold with complex dimension $n$ and $F\>X$ a
   Hermitian vector bundle with Chern connection $\nabla=\nabla'+\nabla''$. At first, by Hodge theory on vector bundles( Lemma
  \ref{hodge}), we observe that if $\alpha\in \Om^{n,0}(X,F)$, and it
  has no harmonic part, then $v=\nabla'^*\G'\alpha$ is a solution to
  $ \nabla' v=\alpha$ where $\G'$ is the Green's operator with respect to $\nabla'$.
Moreover, $\nabla''v$ is a primitive $(n-1,1)$ form. We can apply
this observation to  the K\"ahler fibration $p:\sX\>S$. Let
$(t^1,\cdots, t^m)$ be local holomorphic coordinates on the base $S$
centered at some point $s\in S$. Let $\nabla^\sE=\nabla'+\nabla''$
be the Chern connection of the Hermitian vector bundle $(\sE,
h^{\sE})$ over $\sX$ and $\nabla_X=\nabla'_X+\nabla''_X$ be the
restriction of $\nabla^\sE$ on the fiber $\sE_s\>X_s$. For any
\emph{ local holomorphic} section $u$ of $E=p_*(K_{\sX/S}\ts \sE)$,
by the identification stated above, it can be represented by a local
smooth $\sE$-valued $(n,0)$ form  $u$ on $\sX$ with the property
that
$$\nabla' u=dt^i\wedge \nu_i,\ \ \ \nabla''u=dt^j\wedge \eta_j$$
where $\nu_i$ and $\eta_j$ are forms on $\sX$ of bidegree $(n,0)$
and $(n-1,1)$ respectively. It is easy to see that $\nu_i$ and
$\eta_j$ are not uniquely determined as forms on $\sX$, but their
restrictions to fibers are( see Section 2.3 for more details). It is
worth pointing out that, when restricted to each fiber, $[\eta_j]$
is closely related to the Kodaira-Spencer class of the deformation
$\sX\>S$ (See Remark \ref{etaKS}). We set
$$v_i=-\nabla_X'\G'\pi_{\perp}\left(\nu_i\right)$$
where $\pi_\perp=\mathbb I-\pi$ and $\pi:
\Om^{n,0}(X_s,\sE_s)\>H^{n,0}(X_s,\sE_s)$ is the orthogonal
projection on each fiber. At first, we derive a curvature formula
for $E=p_*(K_{\sX/S}\ts \sE)$ by a   simple method (see Theorem
\ref{main000}).

\btheorem\label{main00} Let $\Theta^E$ be the Chern curvature of
$E=p_*(K_{\sX/S}\ts \sE)$. For any local holomorphic section $u$ of
$E$, the curvature $\Theta^E$
 has the following ``negative form": \beq (\sq\Theta^E u,
u)=c_n\int_{X_s}\sq\left\{\Theta^\sE u,u\right\}-(\Delta_X' v_i,
v_j)\cdot (\sq dt^i\wedge d\bar t^j)+(\eta_i, \eta_j)\cdot (\sq
dt^i\wedge d\bar t^j)\label{NN0} \eeq \etheorem

\noindent We shall explain this curvature formula in details in the
following sections, and  also make a simple example in Section $4$
to explain why this ``negative form" is ``natural".

By a decomposition for the second term on the right hand side of
(\ref{NN0}),
$$(\Delta_X' v_i,
v_j)=(\Delta_X' v_i, \Delta'_Xv_j)+(\Delta_X' v_i,
v_j-\Delta_{X}'v_j)$$
 we obtain a curvature form with significant geometric
interpretations and  it is related to deformation theory of vector
bundles.  Let \beq \alpha_i=\Theta^{\sE}\left(\frac{\p}{\p
t^i}\right)\big|_{X_s}\in \Om^{0,1}(X_s, End(\sE_s)). \label{KS1}
\eeq (Note that, if the family is infinitesimally trivial,
$[\alpha_i]\in H^{0,1}(X_s, End(\sE_s))$ is the Kodaira-Spencer
class (\cite[Proposition~1]{Schumacher-Toma92}) of the deformation
$\sE\>\sX\>S$ in the direction of $\frac{\p}{\p t^i}\in T_sS$.) We
observe that $$ \Delta_X' v_i=-\sq\Lambda_g\left(\alpha_i\cup
u\right)$$ when restricted to the fiber $X_s$ where $\Lambda_g$ is
the contraction operator with respect to the K\"ahler metric
$\omega$ on the fiber $X_s$.

\btheorem\label{KSform234} The curvature $\Theta^E$ of
$E=p_*(K_{\sX/S}\ts \sE)$ has the following ``geodesic form":
\begin{eqnarray} (\sq\Theta^E u, u)\nonumber&= &
c_n\int_{X_s}\sq\left\{\Theta^\sE u,u\right\}-(\alpha_i\cup u,
\alpha_j\cup u)\cdot (\sq dt^i\wedge
d\bar t^j)\\
 &&+(\Delta_X' v_i, \Delta'_X v_j-v_j)\cdot (\sq dt^i\wedge d\bar
t^j)\\
\nonumber&&+(\eta_i, \eta_j)\cdot (\sq dt^i\wedge d\bar
t^j)\label{KSform11}.\end{eqnarray} \etheorem

\noindent Rewriting  each line on the right hand side of
 (\ref{KSform11}) a little bit, we reach the following special
case, which is also of particular interest,  since  the first line
in (\ref{KSform11}) is exactly in the geodesic form.
\bcorollary\label{GLcoro1} Let $(\sL, h^{\sL}=e^{-\phi})$ be a
Hermitian line bundle over $\sX$ such that $(\sL|_{X_s},
h^{\sL}_{X_s})$ is positive on each fiber $X_s$. The  curvature
$\Theta^{E_k}$ of  $E_k=p_*(K_{\sX/S}\ts \sL^k)$ has the form:
\begin{eqnarray} (\sq\Theta^{E_k} u, u)\nonumber&= &
c_n\int_{X_s} kc_{i\bar j}(\phi)\{u,u\}(\sq dt^i\wedge
d\bar t^j)\\
 &&+\frac{1}{k}\left( (\Delta_X'+k)^{-1} \left(\nabla''_X\Delta_X' v_i\right),\nabla''_X\Delta_X' v_j \right)\cdot (\sq dt^i\wedge d\bar
t^j)\label{GLA1}\\
\nonumber&&+(\eta_i, \eta_j)\cdot (\sq dt^i\wedge d\bar
t^j),\end{eqnarray}  where $c_{i\bar j}(\phi)$ is given by\beq
c_{i\bar j}(\phi)=\frac{\p^2\phi}{\p t^i\p\bar t^j}-\left\langle
\bp_X\left(\frac{\p\phi}{\p t^i}\right), \bp_X\left(\frac{\p
\phi}{\p t^j}\right)\right\rangle_g\label{Geo1}\eeq \ecorollary

\noindent \bremark \bd\item The curvature formula (\ref{GLA1}) is
derived implicitly in some special cases by different authors (c.f.
\cite{Berndtsson09b}, \cite{Liu-Sun-Yau09}, \cite{Schumacher13}.)
\item In the real parameter case, $$c(\phi)=\ddot\phi-|\bp_X\dot
\phi|^2_g.$$ When $c(\phi)=0$, it is the geodesic equation in the
space of K\"ahler potentials.  For this comprehensive topic, we just
refer the reader to \cite{Semmes92}, \cite{Donaldson99},
\cite{Chen00}, \cite{PhoStu06}, \cite{Berndtsson09b},
\cite{Berndtsson11a} and references therein.

\item  For the vector bundle case, the authors also expect that  the first line on
the right hand side of (\ref{KSform11}), i.e.
$$c_n\int_{X_s}\sq\left\{\Theta^\sE u,u\right\}-(\alpha_i\cup u,
\alpha_j\cup u)\cdot (\sq dt^i\wedge d\bar t^j)$$ can be written
into certain geodesic form in the space of Hermitian metrics on
$\sE$ when $\sE$ has some stability property (see  formula
(\ref{geodl}) for the line bundle case).

\item If $p:\sX\>S$ is the universal curve with genus $g\geq 2$, i.e. $p:\mathcal T_g\>\mathcal
M_g$.  If $\sL=K_{\mathcal T_g/\mathcal M_g}$, one can deduce
Wolpert's curvature formula (\cite{Wolpert86})  for the (dual)
Weil-Petersson metric on $p_*(K^{\ts 2}_{\mathcal T_g/\mathcal
M_g})$ easily from (\ref{Geo1}) (see also \cite{Siu86},
\cite{Liu-Sun-Yau09} \cite{Berndtsson11} and \cite{Schumacher13}).
\item
 When $k=1$, one can use (\ref{GLA1}) to study the
convex and concave property of the logarithm volume functional on a
Fano manifold (\cite{Berndtsson11a}, see also Theorem \ref{convex1},
Proposition \ref{concave}). Intrinsically, it amounts to the
standard $\bp$-estimate $\|\psi\|\leq \|\bp\psi\|$ on functions with
$\int_X\psi=0$ if the Fano manifold is polarized by its
anti-canonical class.

\ed \eremark

As a first application of Theorem \ref{main00}, we obtain

\btheorem\label{main01} Let $\sX\>S$ be infinitesimally trivial. If
there exists a Hermitian metric on $\sE$ which is Nakano-negative
along the base, then $p_*(K_{\sX/S}\ts \sE)$ is
Nakano-negative.\etheorem

Next, we follow Berndtsson's ideas in his remarkable papers
\cite{Berndtsson09a}, \cite{Berndtsson09b}, \cite{Berndtsson11},
\cite{Berndtsson11a} and set
$$\tilde u=u-dt^i\wedge v_i$$
By using  ``Berndtsson's magic formula" $$ c_n\int_{X_s}
\{u,u\}=c_n\int_{X_s} \left\{\tilde u,\tilde u\right\},$$ we obtain
\btheorem\label{po1} The curvature  $\Theta^E$ of $p_*(K_{\sX/S}\ts
\sE)$ has the following ``positive form": \beq (\sq\Theta^E u,
u)=c_n\int_{X_s}\sq \left\{\Theta^\sE\tilde u, \tilde
u\right\}+\left(\nabla''_X v_i+\eta_i, \nabla_X'' v_j+\eta_j\right)
\cdot (\sq dt^i\wedge d\bar t^j)\label{P0} \eeq \etheorem

\noindent When $\sE$ is a line bundle, the curvature formula
(\ref{P0}) is implicitly obtained by Berndtsson in
\cite{Berndtsson09a}, \cite{Berndtsson09b}, \cite{Berndtsson11} and
\cite{Berndtsson11a}. When $\sE$ is a Nakano-positive vector bundle,
a similar formulation seems to be obtained in
\cite{Mourougane-Takayama08} by using Berndtsson's idea, but $v_i$
are not given explicitly. As it is shown, these $v_i$ play a key
role in these curvature formulas and also their applications.

 Let  $c_{i\bar j}$ be
the $\sE$-valued $(n,0)$-form coefficient of $dt^i\wedge d\bar t^j$
in the local expression of
$$\Theta^\sE(u-dt^i\wedge v_i),$$
and  $d_{i\bar j}$ be the $\sE$-valued $(n,0)$-form coefficient of
$dt^i\wedge d\bar t^j$ in the local expression of
$$\sq\nabla''\nabla' \tilde u.$$

 \btheorem The
curvature $\Theta^E$ of $E=p_*(K_{\sX/S}\ts \sE)$ has the following
``compact form": $$ (\sq\Theta^E u,u)=c_n\int_{X_s}\{ d_{i\bar j},
u\}\cdot (\sq dt^i\wedge d\bar t^j).$$ Moreover, if the family
$\sX\>S$ is infinitesimally trivial,
$$ (\sq\Theta^E u,u)=c_n\int_{X_s}\{ d_{i\bar j},
u\}\cdot (\sq dt^i\wedge d\bar t^j)=c_n\int_{X_s}\{ c_{i\bar j},
u\}\cdot (\sq dt^i\wedge d\bar t^j).$$
 \etheorem

\noindent As applications, we  use it to study the degeneracy of the
curvature tensor of  $p_*(K_{\sX/S}\ts \sE)$ under the assumption
that $\sX\>S$ is infinitesimally trivial and $(\sE,h^{\sE})\>\sX$ is
Nakano semi-positive. In this case,  $c_{i\bar j}$ is closely
related to the geometry of the family $\sE\>\sX\>S$. When
$(\sE,h^{\sE}=e^{-\phi})$ is a relatively positive line bundle,
$c_{i\bar j}$ is the same as the geodesic term $c_{i\bar
j}(\phi)(u)$ defined in (\ref{Geo1}) when the curvature degenerates.
Furthermore, when $H^{n,1}(X_s,\sE_s)=0$, we show that $v_i$ are all
holomorphic over the total space $\sX$ and we can use it to
construct holomorphic automorphisms of the family $\sE\>\sX\>S$ and
study the moduli space of projectively flat vector bundles.

We consider a smooth family of projectively flat vector bundles
$(\sE_s, h^{\sE_s})_{s\in S}$ with polarization \beq  \sq
\Theta^{\sE_s}=\omega_g\ts h^{\sE_s}.\label{polar}\eeq
 Let $W_i$ be the dual vector of the Kodaira-Spencer form $\alpha_i$ defined in (\ref{KS1}), i.e. $W_i$ is  an $End(\sE_s)$-valued $(1,0)$ vector field. Then $v_i$, $u$ and the Kodaira-Spencer vectors $W_i$ are related by \beq  i_{W_i}
u=-v_i\eeq when the curvature of $p_*(K_{\sX/S}\ts \sE)$ is
degenerated. In this case, $W_i$ is an $End(\sE_s)$-valued
holomorphic vector field on the fiber. We also see that, the
horizontal lift of $\frac{\p}{\p t^i}$, $$ V_i=\frac{\p}{\p t^i}-W_i
$$ is a (local) $End(\sE)$-valued holomorphic  vector field
over the total space $\sX$. Moreover, the Lie derivatives of the
curvature tensor of $\sE_s$ with respect to $V_i$ are all zero, i.e.
$$\mathscr L_{V_i}\omega_g=0$$   That means, if
the curvature of $p_*(K_{\sX/S}\ts \sE)$ degenerates at some point
$s\in S$, then the family $\sE\>\sX\>S$ moves by an infinitesimal
automorphism of $\sE$ when the base point varies.

  We can formulate it into a  global version. Let $\sX=X\times \mathbb D$, where $\mathbb D$ is a unit
  disk.  Let $\mathbb E_0\>X$ be a holomorphic vector bundle. If  $(\mathbb E_0, h(t))_{t\in \mathbb D}\>X$ is a smooth family of projectively
  flat vector bundles with  polarization (\ref{polar}).  We denote by $\sE$, the pullback family $p^*_2(\mathbb E_0)$ over $p_2:\sX\>X$.

\btheorem
 If the curvature  $\Theta^E$ of  $E=p_*(K_{\sX/\mathbb D}\ts \sE)$ vanishes
 in a small neighborhood of $0\in \mathbb D$, then there
exists a holomorphic vector field $V$ on $X$ with flows $\Phi_t\in
Aut_{H}(X,\mathbb E_0)$ such that $$\Phi_t^*(\omega_{t})=\omega_0$$
for small $t$.
 \etheorem

 \bremark
 We can also use the holomorphic vector field $V$ to study the uniqueness of Hermitian-Einstein metrics on stable bundles, the
stability of the direct image $p_*(\sE)$ and the asymptotic
stability of $p_*(\sE\ts \sL^k)$ for large $k$. We shall carry it
out in the sequel to this paper. \eremark

\vskip 1\baselineskip

\noindent{\bf{Acknowledgement.}} The second named author would like
to thank V. Tosatti, B. Weinkove and S. Zelditch for many helpful
discussions.

\vskip 2\baselineskip

\section{Background materials}
\subsection{Hodge theory on vector bundles}
Let $(E,h)$ be a Hermitian holomorphic vector bundle over the
compact K\"ahler manifold $(X,\omega)$ and $\nabla= \nabla'
+\nabla''$ be the Chern connection on it. Here, we also have the
relation $\nabla''=\bp$. With respect to metrics on $E$ and $X$, we
set
$$\Delta''=\nabla''\nabla''^*+\nabla''^*\nabla'',$$
$$\Delta'=\nabla'\nabla'^*+\nabla'^*\nabla'.$$
Accordingly, we associate the Green operators and harmonic
projections $\mathbb{G}$, $\mathbb{H}$ and $\mathbb{G}'$,
$\mathbb{H}'$ in Hodge decomposition to them, respectively. More
precisely,
$$\mathbb I=\mathbb{H}+\Delta''\circ \mathbb G,~~~~~~~~~~~~\mathbb I=\mathbb{H}'+\Delta'\circ \mathbb G'.$$
 For any $\phi, \psi\in
\Om^{\bullet,\bullet}(X,E)$, there is a \emph{sesquilinear pairing}
\beq \left\{ \phi, \psi\right\}=\phi^\alpha\wedge \bar {\psi^\beta}
\langle e_\alpha, e_\beta\rangle \eeq  if $\phi=\phi^\alpha
e_\alpha$ and $\psi=\psi^\beta e_\beta$ in the local frame
$\{e_\alpha\}$ of $E$. By the metric compatible property, \beq
\p\{\phi,\psi\}=\{\nabla' \phi,\psi\}+(-1)^{p+q}\{\phi,
\nabla''\psi\} \eeq if $\phi\in\Om^{p,q}(X,E)$.

Let $\Theta^E$ be the Chern curvature of $(E,h)$. It is well-known
\beq \Delta''=\Delta'+[\sq\Theta^E,\Lambda_g]\label{BC}\eeq where
$\Lambda_g$ is the contraction operator with respect to the K\"ahler
metric $\omega$. The following observation plays an important role
in our computations.

\blemma\label{hodge}  Let  $E$ be any Hermitian vector bundle over a
compact K\"ahler manifold $(X,\omega)$. For any $\alpha \in
\Om^{n,0}(X,E)$ with no harmonic part with respect to $\Delta''$,
i.e. $\H(\alpha)=0$, then \bd
\item $\H'(\alpha)=0$;
\item The $(n-1,0)$ form $v=\nabla'^*\G' \alpha$ is a solution to the equation
$$
\nabla'v=\alpha;$$ 
\item $\nabla''v$ is primitive.
\ed

\bproof The first statement follows from the Bochner identity on
E-valued $(n,0)$-forms. More precisely, by (\ref{BC})
$$\Delta''\beta=\Delta'\beta$$
for any $\beta \in\Om^{n,0}(X,E)$. Hence $\H'(\beta)=\H(\beta)$.
 For $(2)$, by Hodge decomposition, we have \be
\nabla' v&=&\nabla'\nabla'^{*}\G'(\alpha)\\
&=&\alpha-\H'(\alpha)-\nabla'^{*}\nabla'\G'(\alpha)\\
&=&\alpha-\H'(\alpha)=\alpha. \ee 
 For $(3)$, let $L_g=\omega\wedge $. By Hodge identity $[\nabla'^*,L_g]=-\sq \nabla''$,
\be \omega\wedge \nabla'' v= L_g\nabla''v&=&L_g\nabla''\nabla'^{*}\G'\alpha\\
&=&-L_g\nabla'^{*}\nabla''\G'\alpha\\
&=&\left(-\sq \nabla''-\nabla'^{*}L_g\right)\nabla''\G'\alpha\\
 &=&0
\ee since $L_g\nabla''\G'\alpha$ is an $(n+1,2)$ form. \eproof
\elemma

The following  Rieman-Hodge bilinear  relation will be used
frequently, and the proof of it can be found in \cite[Corollary
~1.2.36]{Huybrechets05} or \cite[Proposition ~6.29]{Voisin02}.
\blemma\label{primitive} If $\phi, \psi \in
\Om^{p,q}(X,E)\subset\Om^n(X,E)$ are primitive, then \beq
(\phi,\psi)=\left(\sq\right)^{n(n-1)+(p-q)}\int_{X}
\{\phi,\psi\}\label{primitive1}\eeq where $(\bullet, \bullet)$ is
the standard inner product (norm) induced by metrics on $X$ and
$E$.\elemma

\vskip 1\baselineskip
\subsection{Positivity of vector bundles}
Let $\{z^i\}_{i=1}^n$ be  the local holomorphic coordinates
  on $X$ and  $\{e_\alpha\}_{\alpha=1}^r$ be a local frame
 of $E$. The curvature tensor $\Theta^E\in \Gamma(X,\Lambda^2T^*X\ts E^*\ts E)$ has the form
 \beq \Theta^E= R_{i\bar j\alpha}^\gamma dz^i\wedge d\bar z^j\ts e^\alpha\ts e_\gamma,\eeq
where $R_{i\bar j\alpha}^\gamma=h^{\gamma\bar\beta}R_{i\bar
j\alpha\bar \beta}$ and \beq R_{i\bar j\alpha\bar\beta}= -\frac{\p^2
h_{\alpha\bar \beta}}{\p z^i\p\bar z^j}+h^{\gamma\bar
\delta}\frac{\p h_{\alpha \bar \delta}}{\p z^i}\frac{\p
h_{\gamma\bar\beta}}{\p \bar z^j}.\eeq Here and henceforth we
 adopt the Einstein convention for summation.

\bdefinition
 A Hermitian vector bundle
$(E,h)$ is said to be \emph{Griffiths-positive}, if for any nonzero
vectors $u=u^i\frac{\p}{\p z^i}$ and $v=v^\alpha e_\alpha$,  \beq
\sum_{i,j,\alpha,\beta}R_{i\bar j\alpha\bar \beta}u^i\bar u^j
v^\alpha\bar v^\beta>0\eeq $(E,h)$ is said to be
\emph{Nakano-positive}, if for any nonzero vector
$u=u^{i\alpha}\frac{\p}{\p z^i}\ts e_\alpha$, \beq
\sum_{i,j,\alpha,\beta}R_{i\bar j\alpha\bar \beta} u^{i\alpha}\bar
u^{j\beta}>0 \eeq $(E,h)$ is said to be \emph{dual-Nakano-positive},
if for any nonzero vector $u=u^{i\alpha}\frac{\p}{\p z^i}\ts
e_\alpha$, \beq \sum_{i,j,\alpha,\beta}R_{i\bar j\alpha\bar \beta}
u^{i\beta}\bar u^{j\alpha}>0 \eeq It is easy to see that $(E,h)$ is
dual-Nakano-positive if and only if $(E^*,h^*)$ is Nakano-negative.
\noindent The notions of semi-positivity, negativity and
semi-negativity can be defined similarly. We say $E$ is
Nakano-positive (resp. Griffiths-positive, dual-Nakano-positive,
$\cdots$), if it admits a Nakano-positive(resp. Griffiths-positive,
dual-Nakano-positive, $\cdots$) metric. \edefinition

\vskip 1\baselineskip

\subsection{Direct image sheaves of vector bundles}
Let $\sX$ be a K\"ahler manifold with dimension $m+n$ and $S$ a
K\"ahler manifold with  dimension $m$. Let $p:\sX\>S$ be a smooth
K\"ahler fibration. That means, for each $s\in S$,
$$X_s:=p^{-1}(\{s\})$$
is a compact K\"ahler manifold with dimension $n$. Let
$(\sE,h^{\sE})\>\sX$ be a Hermitian holomorphic vector bundle. In
the following, we adopt the setting in \cite[Section
~4]{Berndtsson09a}. Consider the space of holomorphic $\sE$-valued
$(n,0)$-forms on $X_s$,
 $$E_s:=H^0(X_s,\sE_s\ts K_{X_s})\cong H^{n,0}(X_s,\sE_s)$$
where $\sE_s=\sE|_{X_s}$. Here, we assume all $E_s$ has the same
dimension. With a natural holomorphic structure,
$$E=\bigcup_{s\in S}\{s\}\times E_s$$ is isomorphic to the  direct
image sheaf $p_*(K_{\sX/S}\ts \sE)$ if $E$ has certain positive
property.

For every point $s\in S$, we can take a local holomorphic coordinate
$(W;t=(t^1,\cdots, t^m))$ centered at $s$ such that $(W; t)$ is a
unit ball in $\mathbb C^m$, and a system of local coordinates
$\mathfrak U=\{(U_\alpha; z_\alpha=(z_\alpha^1,\cdots, z_\alpha^n),
t)\}$ of $p^{-1}(W)\subset \sX$.  We would like to drop the index
$\alpha$ in the sequel when no confusion arises. By the canonical
isomorphism $K_{\sX/S}|_{X_s}\cong K_{X_s}$, we make the following
identification which will be used frequently in the sequel. For more
details, we refer the reader to \cite[Section~4]{Berndtsson09a}  and
\cite[Section ~2]{Mourougane-Takayama08}. \bd
\item a local smooth section $u$ of $E$ over $S$ is an
$\sE$-valued $(n,0)$ form on $X_s$.  In the local holomorphic
coordinates, $(z,t):=(z^1,\cdots, z^n,t^1,\cdots, t^m)$  on $\sX$,
it is equivalent to the fact that $u\wedge dt^1\wedge \cdots \wedge
dt^m$ is a local section of $K_{\sX}$. Hence, for example, if $u'$
is an $\sE$-valued $(n,0)$ form on $\sX$, such that $u'\wedge
dt^1\wedge \cdots \wedge dt^m=0$, then $u+u'$ and $u$ are the same
local smooth section of $E$ over $S$. That means, if we use an
$\sE$-valued $(n,0)$ form $u$ on $\sX$ to represent a given local
smooth section of $E$, in general, $u$ is not unique. Moreover, two
representatives differ by a form $dt^i\wedge \gamma_i$ on $\sX$
where $\gamma_i$ are $(n-1,0)$ forms on $\sX$.

\item $u$ is a local holomorphic section  of $E$ over $S$, if $\bp_\sX
u$ restricted to the zero form on each fiber $X_s$, that is  \beq
\bp_\sX u=\sum_j dt^i\wedge \eta_i\eeq \ed where $\eta_i$ are
$(n-1,1)$ forms when restricted on each fiber $X_s$. Clearly,
$\eta^i$ are not uniquely determined, but their restrictions to
fibers are.

By the identification above, there is a natural metric on the $E$
induced by metrics $h^{\sE_s}$ on $\sE_s$. For any local smooth
section $u$ of $E$, we define a Hermitian metric $h$ on $E$ by

\beq h(u,u)=c_n\int_{X_s} \left\{u,u\right\}\eeq where
$c_n=(\sq)^{n^2}$.

  Next, we want to define the Chern connection for the Hermitian
  holomorphic vector bundle $(E,h)\>S$. Let
  $\nabla^\sE=\nabla'+\nabla''$ be the Chern connection of $(\sE,h^\sE)$
  over the total space $\sX$. Therefore $\nabla''=\bp_\sX$.  For any local smooth section
  $u$ of $E$, it is also an $\sE$-valued  $(n,0)$-form
  on $\sX$. It is obvious that
  \beq \nabla'' u= d\bar t^j\wedge \tau_{\bar j}+dt^i\wedge \eta_i.\eeq
 Similarly, \beq \nabla' u= dt^i\wedge
\nu_i.\eeq  Here, $\tau_{\bar j}$, $\eta_i$ and $\nu_i$ are local
sections over $\sX$, and again, they are not unique on $\sX$, but
there restrictions to fibers are. The following lemma is given in
\cite[Lemma ~4.1]{Berndtsson09a}.

\blemma\label{DefD} Let $D=D'+D''$ be the Chern connection of  the
Hermitian
  holomorphic vector bundle $(E,h)\>S$, then for any local smooth section
  $u$ of $E$,
  \beq D''u=\tau_{\bar j} d\bar t^j,\ \ \ \ D'u=\pi(\nu_i) dt^i\eeq
where $\pi$ is the orthogonal projection \beq \pi:
\Om^{n,0}(X_s,\sE_s)\>H^{n,0}(X_s,\sE_s)\eeq
  \elemma

The following  result is contained in \cite[Lemma
~2.1]{Berndtsson11}. For the sake of completeness, we include a
 proof here.

\blemma For any local holomorphic section $u$ of $E$, we can choose
a representative of $u$, i.e. an $\sE$-valued $(n,0)$ form on $\sX$
such that \beq\nabla'' u=dt^i\wedge \eta_i\eeq where $\eta_i$ are
primitive $(n-1,1)$ forms when restricted on each fiber $X_s$.
\bproof Let $\hat u$ be an arbitrary $(n,0)$ form on $\sX$ which
represents the local \emph{holomorphic} section $u$ of $E$, i.e.
$$\bp_\sX\hat u= dt^i\wedge \hat \eta_i.$$
Let $\omega$ be the $(1,1)$ form on $\sX$ which restricted to the
K\"ahler forms on each fiber $X_s$. Since $\hat u\wedge \omega$ is
an $(n+1,1)$ form on $\sX$, it can be written as
$$\hat u\wedge \omega= dt^i\wedge y_i$$
for some $(n,1)$ forms $y_i$ on $\sX$. Hence
$$ dt^i\wedge \hat \eta_i\wedge \omega=\bp_\sX\hat u\wedge \omega=\bp_\sX(\hat u\wedge \omega)=- dt^i\wedge (\bp_\sX y_i).$$
Then $\hat \eta_i\wedge \omega=-\bp_\sX y_i$ when restricted on each
fiber $X_s$. Let $\hat v_i$ to be any form on $\sX$ such that
$$\hat v_i\wedge  \omega =y_i.$$
If we set $u=\hat u-dt^i\wedge \hat v_i$, then $\bp_\sX u=dt^i\wedge
\eta_i$ where $$\eta_i=\hat \eta_i+\bp_{\sX}\hat v_i.$$ It is
obvious that $\eta_i$ are primitive $(n-1,1)$ forms when restricted
on each fiber $X_s$.
 \eproof \elemma

\bremark \label{etaKS} Let $[k_i]\in H^{0,1}(X_s, T^{1,0}X_s\ts
\sE_s)$ be the Kodaira-Spencer class in the direction of $\p/\p
t^i$. It is shown in \cite[p.543]{Berndtsson09a} and also
\cite[Lemma 2.2]{Berndtsson11}, when restricted to each fiber,
$\eta_i$ and $k_i\cup u$ define the same class in $H^{n-1,1}(X_s,
\sE_s)$. In particular, if $\sX\>S$ is infinitesimally trivial, we
can choose $\eta_i$   to be zero. \eremark

Let $\Theta^E$ be the Chern curvature of $(E,h)\>S$. The following
formula is obvious. \blemma Let $u$ be a local holomorphic section
of $E$ over $S$, then \beq \bp\p(u,u)=(D''D'u,u)-(D'u,D'u)=(\Theta^E
u, u)-(D' u,D'u)\label{com}\eeq \elemma

\vskip 2\baselineskip

To end this section, we list some notations we shall use in the
sequel:\\

$\bullet$ $D=D'+D''$ the Chern connection on the Hermitian vector
bundle $(E,h)\>S$;\\

$\bullet$ $\nabla^\sE=\nabla'+\nabla''$  the Chern connection of $\sE$ over the total space $\sX$; We will also use $\bp_\sX$ for $\nabla''$ if there is no confusion; \\

$\bullet$ To simplify notations, we will denote the Chern connection
$\nabla^{\sE}|_{X_s}$ of the Hermitian vector bundle
$(\sE_s,h^{\sE_s})\>X_s$ by
$\nabla_X=\nabla'_X+\nabla''_X$  although it depends on $s\in S$;\\

$\bullet$ $d=\p+\bp$ the natural decomposition of $d$ on the base
$S$;\\

$\bullet$ $\{\omega_s\}_{s\in S}$ a smooth family of K\"ahler
metrics on $\{X_s\}_{s\in S}$;\\

$\bullet$ $\G'$ the Green's operator for
$\Delta'_X=\nabla'_X\nabla_X'^*+\nabla_X'^*\nabla'_X$;\\

$\bullet$ $\pi: \Om^{n,0}(X_s,\sE_s)\>H^{n,0}(X_s,\sE_s)$ the
orthogonal projection on the fiber;\\

$\bullet$ $\pi_\perp=\mathbb I-\pi$.\\

\vskip 2\baselineskip

\section{Curvature formulas of direct images of vector bundles}

\subsection{A straightforward computation} In this section, we will derive several general curvature formulas for
direct image $E=p_*(K_{\sX/S}\ts \sE)$  by using Lemma \ref{hodge}.

\noindent The following corollary is a special case of Lemma
\ref{hodge}. \bcorollary\label{3.1} For any local section $u$ of $E$
with $\nabla' u=dt^i\wedge \nu_i$, we set \beq
v_i=-\nabla'^*_X\G'\pi_\perp\left(\nu_i\right).\label{vi}\eeq where
$\nu_i$ is restricted on the fiber $X_s$. Then \bd
\item $\nabla'_X v_i=-\pi_\perp\left(\nu_i\right)$;

\item $\nabla''_Xv_i$ is primitive.

 \ed
\ecorollary

\noindent Before computing the curvature tensors of the direct
images, we need a well-known result:

\blemma  $\p$ and $\bp$ commute with the fiber integration. More
precisely,
$$\bp\int_{X_s} \alpha =\int_{X_s} \bp_\sX \alpha,\ \ \ \ \ \ \  \p\int_{X_s} \alpha =\int_{X_s} \p_\sX \alpha$$
for any smooth $\alpha\in\Om^{\bullet,\bullet}(\sX)$. \elemma

\noindent{\textsf{Note that, in this paper, we  make the following
conventions. Let $u$ be a local holomorphic section of
$E=p_*(K_{\sX/S}\ts \sE)$, \bd \item we always choose a
representative, i.e. an $\sE$-valued $(n,0)$ form on $\sX$ such that
$$\nabla'' u=dt^i\wedge \eta_i,\ \ \ \nabla' u=dt^i\wedge \nu_i$$
where $\eta_i$ are primitive $(n-1,1)$ forms when restricted on each
fiber $X_s$, and $\nu_i$ are $(n,0)$ forms when restricted on each
fiber $X_s$.
\item $v_i$ is fixed to be
$-\nabla'^*_X\G'\pi_\perp\left(\nu_i\right)$, and we do not change
it anymore.  \ed}}

\btheorem\label{main000} Let $\Theta^E$ be the Chern curvature of
$E=p_*(K_{\sX/S}\ts \sE)$. For any local holomorphic section $u$ of
$E$, the curvature $\Theta^E$ has the following ``negative form":
\beq (\sq\Theta^E u, u)=c_n\int_{X_s}\sq\left\{\Theta^\sE
u,u\right\}-(\Delta_X' v_i, v_j)\cdot (\sq dt^i\wedge d\bar
t^j)+(\eta_i, \eta_j)\cdot (\sq dt^i\wedge d\bar t^j)\label{NN} \eeq
\bproof Since $u$ is a local holomorphic section of
$E=p_*(K_{\sX/S}\ts \sE)$, it can be represented by a local smooth
$(n,0)$ form on $\sX$  with the property $\nabla'' u=dt^i\wedge
\eta_i$ where $\eta_i$ are local $(n-1,1)$ forms over $\sX$.
Moreover, when restricted to each fiber $X_s$, $\eta_i$ are all
primitive. By comparing the top degrees along the fiber direction,
we conclude that \beq \int_{X_s}\nabla'\{\nabla'' u, u\}=\int_{X_s}
\nabla'_X\{\nabla'' u,u\}=0,\label{cal}\eeq where the second
identity follows from Stokes' theorem. It is equivalent to the fact
that \be -c_n\sq \int_{X_s}\{\nabla'\nabla''u,u\}&=&c_n\sq
(-1)^{n+1}\int_{X_s}\{\nabla''u,\nabla'' u\}\\
&=&c_n\sq
(-1)^{n+1}\int_{X_s}\{dt^i\wedge \eta_i,dt^j\wedge \eta_j\}\\
&=&-c_n\int_{X_s}\{\eta_i,\eta_j\}\cdot (\sq dt^i\wedge d\bar t^j)\\
&=&(\eta_i,\eta_j) \cdot (\sq dt^i\wedge d\bar t^j) \ee where the
last identity follows from Lemma \ref{primitive} since $\eta_i,
\eta_j$ are primitive $(n-1,1)$  forms on $X_s$. By taking the
conjugate, we see
$$\overline{-c_n\sq
\int_{X_s}\{\nabla'\nabla''u,u\}}=c_n\sq
\int_{X_s}\{u,\nabla'\nabla''u\}$$ On the other hand,
$(\eta_i,\eta_j) \cdot (\sq dt^i\wedge d\bar t^j)$ is a real $(1,1)$
form, and we obtain \beq -c_n\sq
\int_{X_s}\{\nabla'\nabla''u,u\}=(\eta_i,\eta_j) \cdot (\sq
dt^i\wedge d\bar t^j)=c_n\sq
\int_{X_s}\{u,\nabla'\nabla''u\}.\label{th1}\eeq

\noindent By curvature formula (\ref{com}),
\be(\sq \Theta^E u, u)\nonumber&=&\sq(D'u,D'u)-\sq\p\bp\|u\|^2\\
\nonumber&=&\sq(D'u,D'u)-c_n \int_{X_s} \sq \nabla' \left\{
\nabla''u,
  u\right\}-c_n\int_{X_s} \sq (-1)^n\nabla'\{u,\nabla'u\}\\&\stackrel{(\ref{cal})}{=}&\sq(D'u,D'u)-c_n(-1)^{n}\int_{X_s}\left\{
\sq \nabla'
 u, \nabla'  u\right\}-c_n\sq\int_{X_s}\{u,\nabla''\nabla'
u\}\\&=& \sq(D'u,D'u)-c_n(-1)^{n}\int_{X_s}\left\{ \sq \nabla'
 u, \nabla'  u\right\}\\&&+c_n\int_{X_s}\sq\left\{\Theta^\sE
u,u\right\}+(\eta_i,\eta_j) \cdot (\sq dt^i\wedge d\bar t^j)\ee
where the last identity follows from (\ref{th1}) and
$\Theta^\sE=\nabla'\nabla''+\nabla''\nabla'$. By definition (Lemma
\ref{DefD}), we have $D'u=\pi(\nu_i)\wedge dt^i$. From the
orthogonal decomposition,
$\nu_i=\pi(\nu_i)+\pi_{\perp}\left(\nu_i\right)$, and the fact
$\nabla' u=dt^i\wedge \nu_i$,
 we see that
\be &&\sq(D'u,D'u)-c_n(-1)^{n}\int_{X_s}\left\{ \sq \nabla'
 u, \nabla'  u\right\}\\&=&\sq(D'u,D'u)-c_n(-1)^{n}\int_{X_s}\left\{ \sq dt^i\wedge \nu_i, dt^j\wedge \nu_j\right\}\\
 &=&-c_n(-1)^{n}\int_{X_s}\left\{ \sq dt^i\wedge \pi_\perp(\nu_i), dt^j\wedge \pi_\perp(\nu_j)\right\}\\&=&-c_n\int_{X_s}\{\nabla_X' v_i,\nabla_X' v_j\} (\sq dt^i\wedge d\bar
 t^j).\ee
where we use the fact that $-\pi_{\perp}(\nu_i)=\nabla_X' v_i$ in
Corollary \ref{3.1}. Since $\nabla'_X v_i$ are top $(n,0)$ forms on
each fiber, and so primitive. On the other hand, by formula
(\ref{vi}), we know $\nabla'^*_X v_i=0$. Therefore, by Riemann-Hodge
bilinear relation (\ref{primitive1}), we obtain \be
\sq(D'u,D'u)-c_n(-1)^{n}\int_{X_s}\left\{ \sq \nabla'
 u, \nabla'  u\right\}&=&-(\nabla_X' v_i, \nabla_X'v_j)\cdot (\sq dt^i\wedge d\bar
t^j)\\&=&-(\Delta_X' v_i, v_j)\cdot (\sq dt^i\wedge d\bar t^j).\ee
Now the curvature formula (\ref{NN}) follows. \eproof \etheorem


\noindent As a straightforward consequence of Theorem \ref{main000},
we obtain \bcorollary\label{Gcoro} The curvature $\Theta^E$ has the
following form:
\begin{eqnarray} (\sq\Theta^E u, u)\nonumber&= &
c_n\int_{X_s}\sq\left\{\Theta^\sE u,u\right\}-(\Delta_X' v_i,
\Delta'_X v_j)\cdot (\sq dt^i\wedge
d\bar t^j)\\
 &&+(\Delta_X' v_i, \Delta'_X v_j-v_j)\cdot (\sq dt^i\wedge d\bar
t^j)+(\eta_i, \eta_j)\cdot (\sq dt^i\wedge d\bar
t^j)\label{G}.\end{eqnarray} \ecorollary

In the following, we want to interpret the second term on the right
hand side of (\ref{G}) into a geometric quantity. Let \beq
\alpha_i=\Theta^{\sE}\left(\frac{\p}{\p t^i}\right)\big|_{X_s}\in
\Om^{0,1}(X_s, End(\sE_s)). \label{KS2} \eeq  If the family is
infinitesimally trivial,  $[\alpha_i]\in H^{0,1}(X_s, End(\sE_s))$
is the Kodaira-Spencer class
(\cite[Proposition~1]{Schumacher-Toma92}) of the deformation
$\sE\>\sX\>S$ in the direction of $\frac{\p}{\p t^i}\in T_sS$. We
need to point out that, all computations are restricted to the fixed
fiber $X_s(=X)$. By Hodge identity $[\Lambda_g, \nabla''_X]=-\sq
\nabla'^*_X$, $\nabla'^*_Xv_i=0$ and $\nabla''_X\left(
\pi(\nu_i)\right)=0$, we get
\begin{eqnarray} \Delta'_X v_i\nonumber&=&\nabla'^*_X\nabla'_X
v_i=\sq \Lambda_g\nabla''_X\nabla'_X v_i\\ (\text{Corollary
\ref{3.1}}) &=&-\sq \Lambda_g \nabla''_X \pi_\perp(\nu_i)\\
\nonumber&=&-\sq \Lambda_g\nabla''_X\nu_i.\end{eqnarray}
%
\noindent On the other hand, by the relation $\nabla' u=dt^i\wedge
\nu_i$, it is obvious that when restricted to $X_s$,
$$\nu_i=\nabla'_i u,$$ where we adopt the notation that $\nabla'_i
u:=(\nabla' u)(\frac{\p}{\p t^i})$. We get \begin{eqnarray} \Delta_X
v_i= -\sq \Lambda_g\nabla''_X\nu_i\nonumber&=&-\sq
\Lambda_g\left(\nabla''_X\nabla'_i+\nabla'_i\nabla''_X\right)
u+\sq\Lambda_g \nabla'_i\nabla''_X u\\&=&-\sq\Lambda_g(\alpha_i\cup
u)+\sq\Lambda_g \nabla'_i\nabla''_X u\label{KSequ}\\
\nonumber&=&-\sq\Lambda_g(\alpha_i\cup u).
\end{eqnarray} where the last step follows from the
fact that $\sq\Lambda_g \nabla'_i\nabla''_X u$ is zero when
restricted to $X_s$. In fact,$$ \nabla''_X u=\nabla''u-d\bar
t^\ell\wedge \frac{\p u}{\p\bar t^\ell}=dt^j\wedge \eta_j-d\bar
t^\ell\wedge \frac{\p u}{\p\bar t^\ell}.$$

By Corollary \ref{Gcoro}, we obtain the following :
\btheorem\label{KSform123} The curvature $\Theta^E$ has the
following ``geodesic form":
\begin{eqnarray} (\sq\Theta^E u, u)\nonumber&= &
c_n\int_{X_s}\sq\left\{\Theta^\sE u,u\right\}-(\alpha_i\cup u,
\alpha_j\cup u)\cdot (\sq dt^i\wedge
d\bar t^j)\\
 &&+(\Delta_X' v_i, \Delta'_X v_j-v_j)\cdot (\sq dt^i\wedge d\bar
t^j)+(\eta_i, \eta_j)\cdot (\sq dt^i\wedge d\bar
t^j)\label{KSform}.\end{eqnarray} \etheorem

\noindent Next we want to explain why (\ref{KSform}) is called a
``geodesic form" by a little bit more computations in some special
cases. Let $(\sE, e^{-\phi})$ be a relative positive line bundle
over $\sX$, and we
 set $\omega_g=\sq \p_X\bp_X \phi$ on each fiber.

\bcorollary\label{GLcoro} Let $(\sL, h^{\sL}=e^{-\phi})$ be a
Hermitian line bundle over $\sX$ such that $(\sL|_{X_s},
h^{\sL}_{X_s})$ is positive on each fiber $X_s$. Then the curvature
$\Theta^E$ of  $E=p_*(K_{\sX/S}\ts \sL)$ has the   form:
\begin{eqnarray} (\sq\Theta^E u, u)\nonumber&= &
c_n \int_{X_s} c_{i\bar j}
(\phi)\left\{u,u\right\}\cdot (\sq dt^i\wedge d\bar t^j)\\
 &&+\left(
(\Delta_X'+1)^{-1}\left( \nabla''_X\Delta_X' v_i\right),
\nabla''_X\Delta_X' v_j\right)\cdot (\sq dt^i\wedge d\bar t^j)\\
\nonumber&&+(\eta_i, \eta_j)\cdot (\sq dt^i\wedge d\bar
t^j)\label{GL}.\end{eqnarray} where $c_{i\bar j}(\phi)$ is given
by\beq c_{i\bar j}(\phi)=\frac{\p^2\phi}{\p t^i\p\bar
t^j}-\left\langle \bp_X\left(\frac{\p\phi}{\p t^i}\right),
\bp_X\left(\frac{\p \phi}{\p
t^j}\right)\right\rangle_g\label{Geo}\eeq \bproof By formula
(\ref{KSequ}), we obtain \beq c_n\int_{X_s}\sq\left\{\Theta^\sL
u,u\right\}-(\Delta_X' v_i, \Delta'_X v_j)\cdot (\sq dt^i\wedge
d\bar t^j)=c_n \int_{X_s} c_{i\bar j} (\phi)\left\{u,u\right\}\cdot
(\sq dt^i\wedge d\bar t^j) .\label{geodl}\eeq
 On the other hand, we see that, the
second line on the right hand side of (\ref{G}) is non-negative. In
fact, by formula (\ref{BC}) and the fact $\omega_g=\sq
\Theta^{\sL}|_{X_s}$, we have $\Delta''_Xv_i=\Delta'_X v_i-v_i$
since $v_i$ are $(n-1,0)$ forms. Moreover, $\nabla''_X\Delta''_X
v_i+\nabla''_X v_i=\nabla''_X\Delta'_X v_i$. Therefore, \be
(\Delta_X' v_i, \Delta'_X v_j-v_j)&=&(\Delta''_X v_i+v_i, \Delta''_X
v_j)=(\Delta''_X v_i, \Delta''_X v_j)+(v_i, \Delta''_X
v_j)\\
&=&(\Delta''_X v_i, \Delta''_X v_j)+(\Delta''_X  v_i, v_j)=(\Delta''
v_i,\Delta' v_j)\\
&=&(\nabla'' v_i, \nabla''\Delta_X' v_j).\ee  Similarly, by
(\ref{BC}), we have $\Delta'_X(\nabla''_Xv_i)=\Delta''_X(\nabla''_X
v_i)$ since $\nabla''_X v_i$ are $(n-1,1)$ forms on the fiber.
Therefore, $$\Delta'_X(\nabla''_Xv_i)=\Delta''_X(\nabla''_X
v_i)=\nabla''(\Delta''_X v_i)=\nabla''_X(\Delta'_X v_i -v_i)$$ which
is equivalent to
$$(\Delta_X'+1)(\nabla''_X v_i)= \nabla''_X\Delta_X' v_i,$$ or equivalently, $$\nabla''_X v_i=(\Delta_X'+1)^{-1}(\nabla''_X\Delta_X'
v_i).$$ Hence, we obtain
$$(\Delta_X' v_i, \Delta'_X v_j-v_j)=(\nabla'' v_i, \nabla''\Delta_X' v_j)=\left(
(\Delta_X'+1)^{-1}\left( \nabla''_X\Delta_X' v_i\right),
\nabla''_X\Delta_X' v_j\right).$$ \eproof \ecorollary

\bremark Similar formulas are also obtained in \cite{Liu-Sun-Yau09},
\cite{Berndtsson11} and \cite{Schumacher13}. \eremark

Similarly, we  get the ``quantization" version:

\bproposition The  curvature $\Theta^{E_k}$ of $E_k=p_*(K_{\sX/S}\ts
\sL^k)$ has the following form:
\begin{eqnarray} (\sq\Theta^{E_k} u, u)\nonumber&= &
c_n\int_{X_s} kc_{i\bar j}(\phi)\{u,u\}(\sq dt^i\wedge
d\bar t^j)\\
 &&+\frac{1}{k}\left( (\Delta_X'+k)^{-1} \left(\nabla''_X\Delta_X' v_i\right),\nabla''_X\Delta_X' v_j \right)\cdot (\sq dt^i\wedge d\bar
t^j)\\
\nonumber&&+(\eta_i, \eta_j)\cdot (\sq dt^i\wedge d\bar
t^j).\label{GLA}\end{eqnarray} where $c_{i\bar j}(\phi)$ is defined
in (\ref{Geo}). \bproof By Theorem \ref{main000}, we rewrite the
curvature formula as
\begin{eqnarray} (\sq\Theta^{E_k} u, u)\nonumber&= &
c_n\int_{X_s}\sq\left\{\Theta^{\sL^k}
u,u\right\}-\frac{1}{k}(\Delta_X' v_i, \Delta'_X v_j)\cdot (\sq
dt^i\wedge
d\bar t^j)\\
 &&+\frac{1}{k}(\Delta_X' v_i, \Delta'_X v_j-kv_j)\cdot (\sq dt^i\wedge d\bar
t^j)\\
\nonumber&&+(\eta_i, \eta_j)\cdot (\sq dt^i\wedge d\bar
t^j).\end{eqnarray} As similar as the arguments in Corollary
\ref{GLcoro}, we deduce
$$c_n\int_{X_s}\sq\left\{\Theta^{\sL^k}
u,u\right\}-\frac{1}{k}(\Delta_X' v_i, \Delta'_X v_j)\cdot (\sq
dt^i\wedge d\bar t^j)=c_n\int_{X_s} kc_{i\bar j}(\phi)\{u,u\}(\sq
dt^i\wedge d\bar t^j)$$ and
$$(\Delta_X' v_i, \Delta'_X v_j-kv_j)=\left( (\Delta_X'+k)^{-1} \left(\nabla''_X\Delta_X' v_i\right),\nabla''_X\Delta_X' v_j \right).$$
Hence (\ref{GLA}) follows.
 \eproof \eproposition

\vskip 1\baselineskip

\subsection{Computations by using Berndtsson's magic formula} In this subsection, we will derive several curvature formulas for
direct image sheaf $p_*(K_{\sX/S}\ts \sE)$ following the ideas in
\cite{Berndtsson09a}, \cite{Berndtsson09b}, \cite{Berndtsson11} and
\cite{Berndtsson11a}. However, we do not make any assumption on the
curvature of $\sE$. We only make use of the following ``Berndtsson's
magic formula":

\blemma Let $u$ be a local {\bf{smooth}} section of $E$. If  $\tilde
u=u-dt^i\wedge v_i$,  then \beq c_n\int_{X_s} \{u,u\}=c_n\int_{X_s}
\left\{\tilde u,\tilde u\right\}\label{magic}\eeq\bproof It follows
by comparing the $(n,n)$-forms along the fiber $X_s$.\eproof \elemma

\noindent In the following, $u$ shall be a local \emph{holomorphic}
section to $E$, i.e. $\nabla''u=dt^i\wedge \eta_i$. Moreover,  we
set \beq \tilde u=u-dt^i\wedge v_i\eeq and thus fixed. Recall that,
$v_i=-\nabla'^*_X\G'\pi_\perp\left(\nu_i\right)$ as defined in
(\ref{vi}). It is easy to see, $\nabla'' v_i=\nabla''_X v_i+d\bar
t^j\wedge \frac{\p v_i}{\p \bar t^j}$, and so\beq \nabla''\tilde
u=\nabla''u+dt^i\wedge \nabla'' v_i=dt^i\wedge (\eta_i+ \nabla''_X
v_i)+dt^i\wedge d\bar t^j\wedge \frac{\p v_i}{\p \bar
t^j}\label{(0,1)} \eeq and similarly,  \be \nabla' \tilde
u&=&\nabla'u+dt^i\wedge \nabla' v_i\\&=& dt^i\wedge \nu_i+dt^i\wedge
\nabla'v_i\\&=&dt^i\wedge \left(\nu_i+\nabla'_X
v_i\right)+dt^i\wedge dt^k\wedge \nabla'_k v_i\\&=&dt^i\wedge
\pi(\nu_i)+dt^i\wedge dt^k\wedge \nabla'_k v_i\ee since $\nabla'_X
v_i=-\pi_\perp(\nu_i)$. To make the above formula into a compact
form, we define \beq \mu_i:=\pi(\nu_i)\eeq and so \beq\nabla' \tilde
u=dt^i\wedge \mu_i+dt^i\wedge dt^k\wedge \nabla'_k
v_i\label{(1,0)}\eeq
 Therefore,
\begin{eqnarray} \nabla''\nabla'\tilde
u\nonumber&=&\nabla''(dt^i\wedge \mu_i+dt^i\wedge dt^k\wedge
\nabla'_k v_i)\\&=&-dt^i\wedge d\bar t^j\wedge \frac{\p\mu_i}{\p\bar
t^j}+dt^i\wedge dt^k\wedge \nabla''\nabla'_k
v_i\label{(0,1)(1,0)}.\end{eqnarray}

\noindent By curvature formula (\ref{com}) and the magic formula
(\ref{magic}), we obtain
\begin{eqnarray} (\sq \Theta^E u, u)\nonumber&=&\sq(D'u,D'u)-\sq\p\bp\|u\|^2\\ \nonumber&=&\sq(D'u,D'u)-\sq\p\bp\|\tilde u\|^2\\
\nonumber&=&\sq(D'u,D'u)-c_n(-1)^{n}\int_{X_s}\left\{ \sq \nabla'
\tilde u, \nabla' \tilde u\right\}\\&&-c_n \int_{X_s}\sq\left\{
\nabla'\nabla''\tilde u ,\tilde u\right\}-c_n(-1)^{n+1}\int_{X_s}
\sq
\left\{\nabla'' \tilde u, \nabla'' \tilde u\right\}\label{curvature0}\\
\nonumber&&-c_n \int_{X_s} \sq \left\{\tilde u, \nabla''\nabla'
\tilde u\right\}
\end{eqnarray}

{\bf{Claim. }} The first line and second line on the right hand side
of (\ref{curvature0}) are all zero, i.e. \beq
\sq(D'u,D'u)-c_n(-1)^{n}\int_{X_s}\left\{ \sq \nabla' \tilde u,
\nabla' \tilde u\right\}=0\label{zero1}\eeq and \beq-c_n
\int_{X_s}\sq\left\{ \nabla'\nabla''\tilde u ,\tilde
u\right\}-c_n(-1)^{n+1}\int_{X_s} \sq \left\{\nabla'' \tilde u,
\nabla'' \tilde u\right\}=0 \label{zero2}\eeq

\bproof In fact, thanks to (\ref{(1,0)}), we have \be
-c_n(-1)^{n}\int_{X_s}\left\{ \sq \nabla' \tilde u, \nabla' \tilde
u\right\}&=&-c_n(-1)^n\int_{X_s} \sq \left\{dt^i\wedge \mu_i, d
t^j\wedge
 \mu_j \right\}\\&=&-c_n\int_{X_s} \left\{\mu_i, \mu_j\right\}\cdot
(\sq dt^i\wedge d\bar t^j)\ee On the other hand,
$$D'u=\pi(\nu_i)\wedge dt^i=\mu_i\wedge dt^i.$$
Hence $$\sq (D'u,D'u)=c_n\int_{X_s} \left\{\mu_i, \mu_j\right\}\cdot
(\sq dt^i\wedge d\bar t^j).$$ We complete the proof of
(\ref{zero1}). On the other hand,
$$-c_n
\int_{X_s}\sq\left\{ \nabla'\nabla''\tilde u ,\tilde
u\right\}-c_n(-1)^{n+1}\int_{X_s} \sq \left\{\nabla'' \tilde u,
\nabla'' \tilde u\right\}=-c_n\int_{X_s}\nabla'\{\nabla''\tilde
u,\tilde u\}.$$ By formula (\ref{(0,1)}), $\{\nabla''\tilde u,\tilde
u\}$ is an $(n,n+1)$ form on the total space $\sX$, and contains
factors $dt^i$. To make a volume form on the fiber $X_s$, we obtain
$$-c_n\int_{X_s}\nabla'\{\nabla''\tilde
u,\tilde u\}=-c_n\int_{X_s}\nabla'_X\{\nabla''\tilde u,\tilde
u\}=0,$$ by Stokes' theorem.
 Hence, (\ref{zero2}) follows.
\eproof

By taking conjugate of the real  forms,  the curvature formula
(\ref{curvature0}) can be written as
\begin{eqnarray} (\sq \Theta^E u, u)\nonumber&=&-c_n \int_{X_s} \sq \left\{\tilde u,
\nabla''\nabla' \tilde u\right\}\\
\nonumber&=&c_n\int_{X_s}\sq \left\{\nabla''\nabla'\tilde u, \tilde u\right\}\\
&=&\label{curvature1}c_n\int_{X_s}\sq \left\{\Theta^\sE\tilde u,
\tilde u\right\}-c_n\int_{X_s}\sq \left\{\nabla'\nabla''\tilde u,
\tilde u\right\}\\
\nonumber&\stackrel{(\ref{zero2})}{=}&c_n\int_{X_s}\sq
\left\{\Theta^\sE\tilde u, \tilde
u\right\}+c_n(-1)^{n+1}\int_{X_s}\sq
\left\{\nabla''\tilde u, \nabla''\tilde u\right\}\\
\nonumber&=&c_n\int_{X_s}\sq \left\{\Theta^\sE\tilde u, \tilde
u\right\}-c_n \int_{X_s} \left\{\eta_i+\nabla''_X v_i,
\eta_j+\nabla_X'' v_j\right\} \cdot (\sq dt^i\wedge d\bar t^j)
\end{eqnarray}
where the last identity follows from formula (\ref{(0,1)}). Note
that, since $\nabla''_X v_i$ and $\eta_i$ are primitive $(n-1,1)$
forms on $X_s$, by Riemann-Hodge bilinear relation (Lemma
\ref{primitive}), \beq -c_n \int_{X_s} \left\{\eta_i+\nabla''_X v_i,
\eta_j+\nabla_X'' v_j\right\} =(\eta_i+\nabla''_X
v_i,\eta_j+\nabla''_X v_j).\eeq  In summary, we obtain,
\btheorem\label{po} The curvature $\Theta^E$ of
$p_*(K_{\sX/S}\ts\sE)$ has the following ``positive form": \beq
(\sq\Theta^E u, u)=c_n\int_{X_s}\sq \left\{\Theta^\sE\tilde u,
\tilde u\right\}+\left(\eta_i+\nabla''_X v_i, \eta_j+\nabla_X''
v_j\right) \cdot (\sq dt^i\wedge d\bar t^j)\label{P} \eeq \etheorem

Let  $c_{i\bar j}$ be the $\sE$-valued $(n,0)$-form coefficient of
$dt^i\wedge d\bar t^j$ in the local expression of
$$\sq\Theta^\sE(\tilde u)$$
and  $d_{i\bar j}$ be the $\sE$-valued $(n,0)$-form coefficient of
$dt^i\wedge d\bar t^j$ in the local expression of
$$\sq\nabla''\nabla' \tilde u.$$

\btheorem\label{compact} The curvature $\Theta^E$ of
$p_*(K_{\sX/S}\ts \sE)$ has the following ``compact form": \beq
(\sq\Theta^E u,u)=c_n\int_{X_s}\sq \{\nabla''\nabla' \tilde u,
u\}=c_n\int_{X_s}\{ d_{i\bar j}, u\}\cdot (\sq dt^i\wedge d\bar
t^j).\label{compctform}\eeq

 \bproof By formula (\ref{curvature1}),
 $$(\sq\Theta^E u,u)=-c_n\int_{X_s}\sq \{\tilde u, \nabla''\nabla' \tilde u\}=c_n\int_{X_s}\sq \{\nabla''\nabla' \tilde u, \tilde u\}.$$
We obtain from formula (\ref{(0,1)(1,0)}) that
$$\int_{X_s} \{\nabla''\nabla'\tilde u, dt^i\wedge v_i\}=0$$
for degree reasons. Therefore $$(\sq\Theta^E u,u)=c_n\int_{X_s}\sq
\{\nabla''\nabla' \tilde u, u\}.$$ By degree reasons again, we
obtain the last identity in (\ref{compctform}).
 \eproof
\etheorem

\bremark If $\sX\>S$ is infinitesimally trivial, we can choose a
representative $u$ such that $\nabla''u=0$ (i.e. $\eta_i=0$, see
Remark \ref{etaKS}) and so \be \int_{X_s}
\{\nabla'\nabla'' \tilde u, u\}&=&-\int_{X_s} \{\nabla'(dt^i\wedge v_i), u\}\\&=&-\int_{X_s}\{\nabla_X'(dt^i\wedge v_i), u\}\\
&=&(-1)^n\int_X\{dt^i\wedge v_i, \nabla''_Xu\}=0. \ee Therefore
 \beq
(\sq\Theta^E u,u)=c_n\int_{X_s}\sq \{\Theta^{\sE} \tilde u,
u\}=c_n\int_{X_s}\{ c_{i\bar j}, u\}\cdot (\sq dt^i\wedge d\bar
t^j)\label{compctform2}.\eeq \eremark

\vskip 2\baselineskip

\section{Curvature positivity and negativity for direct images of vector bundles}

As applications of curvature formulas derived in Section $3$, at
first, we obtain \btheorem\label{main1} Let $\sX\>S$ be an
infinitesimally trivial proper holomorphic fibration. If there
exists a Hermitian metric on $\sE$ which is Nakano-negative along
the base, then $p_*(K_{\sX/S}\ts \sE)$ is Nakano-negative. \bproof
It follows from Theorem \ref{main000}. Here, the curvature formula
(\ref{NN}) is reduced to \beq  (\sq\Theta^E u,
u)=c_n\int_{X_s}\sq\left\{\Theta^\sE u,u\right\}-(\Delta_X' v_i,
v_j)\cdot (\sq dt^i\wedge d\bar t^j)\label{NNNN} \eeq We set
$u=u^{i\alpha} \frac{\p}{\p t^i} \ts e_\alpha$. Naturally,
$e_\alpha$ can be viewed as a local holomorphic section of
$H^{n,0}(X_s,\sE_s)$. We set
$$\Theta^{\sE}_{i\bar j\alpha\bar\beta}=\Theta^{\sE}\left(\frac{\p}{\p t^i} \ts e_\alpha, \frac{\p}{\p t^j} \ts e_\beta\right),$$
 then by (\ref{NNNN}), \beq
\Theta^{E}_{i\bar j\alpha\bar\beta}
u^{i\alpha}\bar{u}^{j\beta}=c_n\int_{X_s}\Theta^{\sE}_{i\bar j
\alpha\bar\beta}\{u^{i\alpha},u^{j\beta}\}-(\Delta'_X v,
v)\label{NNN}\eeq where
$v=-\sum_i\nabla'^*_X\G'\pi_\perp\left(\nu_i^{i\alpha}\ts
e_\alpha\right)$. If $\sE$ admits a Hermitian metric
 which is Nakano-negative along the base, then the first term
 in the formula (\ref{NNN}) is negative.
 \eproof\etheorem

\bcorollary[\cite{LSYang}] If $(E,h)$ is a Griffiths-positive vector
bundle, then $E\ts \det E$ is both Nakano positive and dual
Nakano-positive. \bproof The Nakano-positivity is
well-known(\cite{Demailly}, \cite{Berndtsson09a}). Now we prove the
dual Nakano-positivity. Let $L=\sO_{\P(E)}(1)$ be the tautological
line bundle of $\P(E)$. Note that $\sO_{\P(E^*)}(1)$ is ample, but
$L$ is not. The metric on $(E,h)$ induces a metric on $L$ which is
negative along the base(\cite[Chapter V, formula 15.15]{Demailly},
\cite[(2.12)]{LSYang}).
 On the other hand, it is easy to see
$$E^*\ts \det E^*=p_*( K_{\P(E)/S}\ts L^{r+1})$$
where $p:\P(E)\>S$ is the projection. Hence, by Theorem \ref{main1},
$E^*\ts \det E^*$ is Nakano-negative, or equivalently, $E\ts \det E$
is dual Nakano-positive. \eproof \ecorollary

\noindent Similarly, for the Nakano-positivity, it follows from
Theorem \ref{po} and the proof is similar to that of Theorem
\ref{main1}.

\bcorollary[\cite{Mourougane-Takayama08}]\label{MT} \label{}
$p_*(K_{\sX/S}\ts \sE)$ is Nakano-positive if $\sE$ is
Nakano-positive. \ecorollary

\bcorollary[\cite{Berndtsson09a}]\label{Be} Let $\sL$ be a line
bundle over $\sX $. Then $p_*(K_{\sX/S}\ts \sL)$ is Nakano-positive
if $\sL$ is ample. \ecorollary

\bremark Note that,  in Corollary \ref{MT} and Corollary \ref{Be},
the family $\sX\>S$ are {\bf{not}} necessarily infinitesimally
trivial since the term related to the Kodaira-Spencer class
$$\left(\eta_i+\nabla''_X v_i, \eta_j+\nabla_X''
v_j\right) \cdot (\sq dt^i\wedge d\bar t^j)$$
 is nonnegative.

\eremark

Let $X$ be a compact Fano manifold and  $h(t)=e^{-\phi(t)}$ be a
family of positive metrics on $L=-K_X$. Let $(z^1,\cdots, z^n)$ be
the local holomorphic coordinates on $X$. We set the \emph{local}
volume form $$dV_{\mathbb C}=(\sq)^n(dz^1\wedge d\bar
z^1+\cdots+dz^n\wedge d\bar z^n)^n.$$ It is easy to see that
$$e^{-\phi} dV_{\mathbb C}$$
is a family of globally defined volume forms of $X$. Berndtsson in
\cite{Berndtsson11a} considers the logarithm volume\beq \mathcal
F(t)=-\log\left(\int_X e^{-\phi} dV_{\mathbb C}\right)\eeq and
deduces that

\btheorem[\cite{Berndtsson11a}]\label{convex1} If $e^{-\phi(t)}$ is
a  subgeodesics in the K\"ahler cone $\mathcal K_{L}$ of the class
$c_1(L)$, i.e.
$$c(\phi)=\ddot\phi-|\bp_X\dot\phi|^2\geq 0,$$ then $\mathcal F(t)$ is
convex.
 \etheorem

\noindent In fact, Theorem \ref{convex1} can be obtained  easily
from Corollary \ref{GLcoro}, following the setting in
\cite{Berndtsson11a}. To formulate it efficiently, we use complex
parameter $t$ in the unit disk $\mathbb D\subset \mathbb C$. When we
consider the direct image bundle $E=p_*(K_X\ts L)$, it is a trivial
line bundle since $L=-K_X$ and $H^0(X,K_X\ts L)\cong \mathbb C$.
Since $E$ is trivial, there is a constant section $u=1 e_E$ of $E$,
and it is identified as  a holomorphic section $u$ of
$H^{n,0}(X,L)$,
$$ u=dz^1\wedge \cdots \wedge dz^n\ts e$$ where $e=\frac{\p }{\p
z^1}\wedge \cdots \wedge \frac{\p}{\p z^n}$. Hence
$$ \|u\|^2=c_n\int_X \{u,u\}=\int_X e^{-\phi} dV_{\mathbb C}$$
On the other hand, it is obvious that \beq \|u\|^2 \sq\p\bp \mathcal
F=(\sq \Theta^{E}u,u)\label{bconv}\eeq Hence, if $c(\phi)\geq 0$, by
Corollary \ref{GLcoro}, $\sq\p\bp\mathcal F$ is Hermitian
semi-positive. In real parameters, it says that $\mathcal F$ is
convex.

As a partial converse to Berndtsson's result, we have \bproposition
\label{concave} Let $e^{-\phi(t)}$ be a curve in the K\"ahler cone
$\mathcal K_{L}$.  If $\phi(t)$ is concave in $t$, then so is
$\mathcal F(t)$. \bproof It follows from Theorem \ref{main000}. In
fact, $\ddot\phi\leq 0$ implies the first term on the right hand
side of (\ref{NN}) is negative. Note that, in this case, the family
is a trivial family and so the third term on the right hand side of
(\ref{NN}) is zero. Therefore, $(\sq \Theta^{E}u,u)\leq 0$. By
formula (\ref{bconv}), we see $\mathcal F$ is superharmonic and in
the real case, it is concave.
 \eproof
 \eproposition

\noindent We can also see how Theorem \ref{convex1} and Proposition
\ref{concave} work  by the following  simple example. At first, we
fix a positive metric $e^{-\phi(0)}$ in $c_1(L)$ and set
$$\phi(t)=f(t)+\phi(0)$$
where $t$ is a real parameter.  It is obvious that $$c(\phi)=\ddot
\phi=\ddot f,\ \ \ \ \ \mathcal F(t)=f(t)+c$$ Hence $\mathcal F$ is
concave if $\phi$ is concave and vice visa.

\indent For the general case, it is not hard to see that both
Theorem \ref{convex1} and Proposition
 \ref{concave} amount to the basic $\bp$-estimate
 \beq \|\dot \phi\|\leq \|\bp_X\dot\phi\|\eeq
if the Fano manifold $X$ is polarized by its anti-canonical class.

\vskip 2\baselineskip

\section{Direct images of projectively flat vector bundles}
In this section we consider an infinitesimal trivial family $\sX\>S$
and assume that the vector bundle $(\sE,h^{\sE})\>\sX$ is Nakano
semi-positive. In this case, we can choose $\eta_i$ to be zero in
all formulas derived in Section $3$.

Let's recall that $c_{i\bar j}$ is the $\sE$-valued $(n,0)$-form
coefficient of $dt^i\wedge d\bar t^j$ in the expression
$$\Theta^\sE(u-dt^i\wedge v_i).$$
There are four (linearly independent) terms in the expression of
$\Theta^\sE(u-dt^i\wedge v_i)$. However, if $\Theta^{\sE}$ is Nakano
semi-positive, then $c_{i\bar j}$ dominates the degeneracy of
$\Theta^\sE(u-dt^i\wedge v_i)$, i.e. $c_{i\bar j}=0$ implies
$\Theta^\sE(u-dt^i\wedge v_i)=0$. This is the content of the next
theorem.

 \btheorem\label{degeneracy} Let $(\sE,h^{\sE})\>\sX$ be Nakano semi-positive.
Then
$$(\sq\Theta^E u,u)=0$$
if and only if $c_{ij}=0$. \bproof Note that if $(\sE,h^{\sE})\>\sX$
is Nakano semi-positive, then by formula (\ref{P}),
$$(\sq\Theta^E u,u)=c_n\int_{X_s}\sq \left\{\Theta^\sE\tilde u,
\tilde u\right\}+\left(\nabla''_X v_i, \nabla_X'' v_j\right) \cdot
(\sq dt^i\wedge d\bar t^j)$$ is a Hermitian semi-positive
$(1,1)$-form. If $(\sq\Theta^E u,u)=0$,  we get
$$c_n\int_{X_s}\sq \left\{\Theta^\sE\tilde u,
\tilde u\right\}=0$$ and so $\Theta^\sE \tilde u=0$. In particular,
$c_{i\bar j}=0$ since $\Theta^\sE \tilde u=\Theta^\sE(u-dt^i\wedge
v_i)$.

 On the other hand,  by Theorem \ref{compact}, if $c_{ij}=0$,
$$(\sq\Theta^E u,u)=c_n\int_{X_s}\{ c_{i\bar j}, u\}\cdot (\sq
dt^i\wedge d\bar t^j)=0.$$
 \eproof \etheorem


Now we continue to analyze the case $(\sq\Theta^E u,u)=0$. In the
following, we use an idea in \cite{Berndtsson11a}.

\blemma If $H^{n,1}(X_s,\sE_s)=0$, then $v_i$ is holomorphic on
$\sX$ for any $i$. \bproof We only need to show $\frac{\p
v_i}{\p\bar t^j}=0$ since $\nabla''_X v_i=0$ is obvious from
curvature formula (\ref{P}) when $(\sq\Theta^E u,u)=0$.

 Next we
claim \beq \nabla'_X \left(\frac{\p v_i}{\p\bar t^j}\right)=c_{i\bar
j}+\frac{\p}{\p\bar t^j}\pi\left(\nu_i\right)\label{holomorphic}\eeq
In fact, \be -\nabla'\nabla''(dt^i\wedge v_i)
&=&-\Theta^\sE(dt^i\wedge v_i)+\nabla''\nabla'(dt^i\wedge v_i)\\
&=&-\Theta^\sE(dt^i\wedge v_i)-\nabla''(dt^i\wedge
\nabla'_Xv_i)-\nabla''(dt^i\wedge dt^k\wedge \nabla'_k v_i)\\
&=&-\Theta^\sE(dt^i\wedge v_i)+\nabla''\left(dt^i\wedge
\left(\nu_i-\pi(\nu_i)\right)\right)-\nabla''(dt^i\wedge
dt^k\wedge \nabla'_k v_i)\\
&=&-\Theta^\sE(dt^i\wedge
v_i)+\nabla''\nabla'u-\nabla''\left(dt^i\wedge
\pi(\nu_i)\right)-\nabla''(dt^i\wedge
dt^k\wedge \nabla'_k v_i)\\
&=&-\Theta^\sE(dt^i\wedge v_i)+\Theta^\sE
u-\nabla'\nabla''u-\nabla''\left(dt^i\wedge
\pi(\nu_i)\right)-\nabla''(dt^i\wedge
dt^k\wedge \nabla'_k v_i)\\
&=&\Theta^\sE(u-dt^i\wedge v_i)-\nabla'(dt^i\wedge
\eta_i)-\nabla''\left(dt^i\wedge
\pi(\nu_i)\right)-\nabla''(dt^i\wedge dt^k\wedge \nabla'_k v_i) \ee
By comparing the coefficients of $dt^i\wedge d\bar t^j$ on both
sides, we get (\ref{holomorphic}). If the curvature is zero, we also
have $c_{i\bar j}=0$. According to different types in Hodge
decomposition, i.e. $\nabla'_X \left(\frac{\p v_i}{\p\bar
t^j}\right)\in Im(\nabla'_X)$ and the  holomorphic $(n,0)$ form
$\frac{\p}{\p\bar t^j}\pi\left(\nu_i\right)\in
Ker(\Delta''_X)=Ker(\Delta'_X)$ (when restricted on $(n,0)$ forms),
 we conclude from (\ref{holomorphic}) that
$$\nabla'_X \left(\frac{\p v_i}{\p\bar t^j}\right)=\frac{\p}{\p\bar t^j}\pi\left(\nu_i\right)=0$$
Therefore, by Hodge relation $[\nabla''^*_X, \omega]=-\sq
\nabla_X'$, we get
 $\nabla''^*_X (\omega\wedge \frac{\p
v_i}{\p\bar t^j})=0$ and $\nabla''_X (\omega\wedge \frac{\p
v_i}{\p\bar t^j})=0$. The cohomology assumption ensures the $(n,1)$
form $\omega\wedge \frac{\p v_i}{\p\bar t^j}=0$ and so $\frac{\p
v_i}{\p\bar t^j}=0$. \eproof \elemma

 In the following, we assume that
$(\sE_s,h^{\sE_s})$ is projectively flat. Hence, the curvature
tensor can be written as (c.f. \cite[p.7]{Kobayashi87}) \beq \sq
\Theta^{\sE_s}=\frac{1}{r}Ric(\det\sE_s)\ts h^{\sE_s}\eeq where $r$
is the rank of $\sE_s$. If $\det \sE_s$ is positive, we set \beq
\omega_g=\frac{1}{r} Ric(\sE_s)=-\frac{\sq}{r} \p_X\bp_X\log
\det(h^{\sE_s}) \eeq as the background K\"ahler metric on each
fiber. Therefore, \beq \sq \Theta^{\sE_s}=\omega_g\ts
h^{\sE_s}\label{pcurv}\eeq

\noindent Recall that $[\alpha_i]\in H^{0,1}(X_s, End(\sE_s))$ is
the Kodaira-Spencer class of the deformation $\sE\>\sX\>S$ in the
direction of $\frac{\p}{\p t^i}$, i.e. \beq
\alpha_i=\Theta^{\sE}\left(\frac{\p}{\p t^i}\right)\big|_{X_s} \in
\Om^{0,1}(X_s, End(\sE_s)). \label{KS} \eeq
 Let $W_i$ be the dual vector of $\alpha_i$, i.e. $W_i$ is an $End(\sE_s)$-valued $(1,0)$-vector field on the fiber $X_s$. Then by
 formulas (\ref{pcurv}) and (\ref{KS}), we have
 \begin{eqnarray} \left(i_{W_i}\omega\right)\wedge u&=&\nonumber\sq \alpha_i\wedge u=\sq \Theta^{\sE}\left(\frac{\p}{\p t_i},u\right)\\&=&\sq \nabla''_X\nu_i\label{defv}
 \end{eqnarray}
since $\nabla_X''u=0$.
 \bproposition We have the relation \beq i_{W_i}
u=-v_i.\label{holov}\eeq Moreover,  $W_i$ is an $End(\sE_s)$-valued
holomorphic vector field on the fiber $X_s$.

 \bproof By formula (\ref{defv}), we obtain \be
\left(i_{W_i}\omega_g\right)\wedge u&=&\sq \nabla''_X\nu_i=-\sq \nabla''_X\nabla'_X v_i\\
&=&-\sq\Theta^{\sE_s}(v_i) \ee since $\nabla''_X v_i=0$. On the
other hand, $\left(i_{W_i}\omega_g\right)\wedge u=\left(i_{W_i}
u\right)\wedge \omega_g$. Hence we obtain (\ref{holov}) by using
(\ref{pcurv}) again. Since $v_i$ and $u$ are all holomorphic on each
fiber, we know $W_i$ is also holomorphic. \eproof \eproposition

We can extend the vector field $\frac{\p}{\p t^i}$ to an
$End(\sE_s)$-valued vector field. We still denote it by
$\frac{\p}{\p t^i}$. Then \beq V_i=\frac{\p}{\p t^i}-W_i \eeq is a
(local) $End(\sE)$-valued holomorphic  vector field over the total
space $\sX$. Let $\mathscr{L}$ be the type $(1,0)$,
$End(\sE)$-valued Lie derivative, then we have \beq \mathscr
L_{V_i}\omega_g=0\label{Lie}\eeq In fact, by relation (\ref{KS}), we
have
$$\mathscr L_{W_i}\omega_g=\nabla'_X (-\sq \alpha)=-\sq \nabla'_X\left(\Theta^\sE\left(\frac{\p}{\p t^i}\right)|_{X_s}\right)=-\nabla_{\frac{\p}{\p t^i}}' \Theta^{\sE_s}$$
Hence, by formula (\ref{pcurv}), we get (\ref{Lie}).  That means, if
the curvature $\Theta^E$ degenerates at some point $s\in S$, then
the family $\sE\>\sX\>S$ moves by an infinitesimal automorphism of
$\sE$.

 We summarize the above into a  global version. Let $\sX=X\times \mathbb D$, where $\mathbb D$ is a unit
  disk.  Let $\mathbb E_0\>X$ be a holomorphic vector bundle. If  $(\mathbb E_0, h(t))_{t\in \mathbb D}\>X$ is a smooth family of projectively
  flat vector bundles.  We assume  $Ric(\det E, h(t))>0$ for
all $t$ and set $\omega_t=-\frac{\sq}{r}\p_X\bp_X\log\det (h(t))$ to
be a smooth family of K\"ahler metrics on $X$. We also denote by
$\sE$, the pullback family $p^*_2(\mathbb E_0)$ over $p_2:\sX\>X$.

\btheorem
 If the curvature  $\Theta^E$ of  $E=p_*(K_{\sX/\mathbb D}\ts \sE)$ vanishes
 in a small neighborhood of $0\in \mathbb D$, then there
exists a holomorphic vector field $V$ on $X$ with flows $\Phi_t\in
Aut_{H}(X,\mathbb E_0)$ such that $$\Phi_t^*(\omega_{t})=\omega_0$$
for small $t$.
 \etheorem

\end{document}